\documentclass[12pt]{article}
\usepackage{amsmath,amsfonts,amstext,amssymb,amsbsy,amsopn,amsthm,eucal}
\usepackage{dsfont}

\begin{document}

\newtheorem*{thm1}{Theorem 1}
\newtheorem*{thm2}{Theorem 2}
\newtheorem*{thm3}{Theorem 3}
\newtheorem*{thm4}{Theorem 4}
\newtheorem*{thm5}{Theorem 5}
\newtheorem*{thm6}{Theorem 6}

\newtheorem*{lemma1}{Lemma 1}
\newtheorem*{lemma2}{Lemma 2}
\newtheorem*{lemma3}{Lemma 3}

\newtheorem*{cor1}{Corollary 1}
\newtheorem*{cor2}{Corollary 2}
\newtheorem*{cor3}{Corollary 3}
\newtheorem*{cor4}{Corollary 4}
\newtheorem*{cor5}{Corollary 5}

\newenvironment{definition}[1][Definition.]{\begin{trivlist}
\item[\hskip \labelsep {\bfseries #1}]}{\end{trivlist}}
\newenvironment{remark}[1][Remark.]{\begin{trivlist}
\item[\hskip \labelsep {\itshape #1}]}{\end{trivlist}}

\title{Some Geometry and Analysis on Ricci Solitons}

\author{Aaron Naber\thanks{Department of Mathematics, Princeton University,
        Princeton, NJ 08540 ({\tt anaber@math.princeton.edu}).}}

\date{\today}
\maketitle
%%% ----------------------------------------------------------------------

\begin{abstract}
The Bakry-\'{E}mery Ricci tensor of a metric-measure space
$(M,g,e^{-f}dv_{g})$ plays an important role in both geometric
measure theory and the study of Hamilton's Ricci flow.  Under a
uniform positivity condition on this tensor and with bounded Ricci
curvature we show the underlying space has finite $f$-volume. As a
consequence such manifolds, including shrinking Ricci solitons, have
finite fundamental group. The analysis can be extended to classify
shrinking solitons under convexity or concavity assumptions on the
measure function.

%The Bakry-\'{E}mery Ricci tensor plays an important role in both
%geometric measure theory and the study of Hamilton's Ricci flow.
%Under various geometric conditions on the Bakry-\'{E}mery Ricci
%tensor we prove topological results about the underlying manifold,
%as well as provide analysis of a natural analogue of the Laplace
%operator on metric measure spaces.
\end{abstract}

%%% ----------------------------------------------------------------------
\section{Introduction}
In this paper we study smooth metric measure spaces $(M,g,f)$, where
$g$ is a smooth complete metric on an $n$ dimensional manifold $M$
and $f$ is a smooth real valued function. We associate to $M$ the
measure $e^{-f}dv_{g}$, where $dv_{g}$ is the Riemannian volume form
on $M$ . The interest of this paper is in studying the
Bakry-\'{E}mery Ricci tensor $Rc_{f}\equiv Rc+\nabla^{2}f$, where
$Rc$ is the usual Ricci tensor and $\nabla^{2}f$ is the hessian of
$f$. We refer the reader to \cite{BE} and to Lott's paper \cite{L}
for more information.

Manifolds with constant Bakry-\'{E}mery tensor have come to be known
as Ricci solitons, and play an important role in the Ricci Flow as
they are the result of certain singularity dilations around finite
time singularities of the Ricci Flow (see \cite{P} and \cite{MT}).
With this in mind we will be interested in studying the following

\begin{definition}
Let $(M,g,f)$ be a smooth metric measure space.  We call $M$ a Ricci
soliton if $Rc_{f}=Rc+\nabla^{2}f=\lambda g$, where
$\lambda\in\mathds{R}$.  We say the soliton is shrinking, steady, or
expanding when $\lambda >0, =0, <0$, respectively.
\end{definition}

Our first result is a form of Myers Theorem for metric measure
spaces with uniform positive lower bounds on the $Rc_{f}$ tensor. It
is well known that such manifolds need not be compact, and in fact
some of the most interesting examples are those which are not. Hence
bounds on the diameter are not reasonable under such a constraint,
however in the following we show that such manifolds do have finite
$f$-volume.

\begin{thm1}
Let $(M,g,f)$ be complete with bounded Ricci curvature. Assume
$Rc_{f} \geq \lambda g$ with $\lambda > 0$.  Then the measure
$e^{-f}dv_{g}$ is finite, and consequently M has finite fundamental
group.
\end{thm1}
\begin{remark}
The finiteness of the fundamental group was proved by Lott in
\cite{L} under the additional assumption that $M$ is compact.
\end{remark}

%The topological behavior of manifolds with positive lower bounds on
%$Rc_{f}$ is not well understood.  The following guarantees for
%noncompact manifolds  with such lower bounds that outside some
%compact region $K\subset\subset M$ the topology of our manifold
%stabilizes.

%\begin{thm2}
%Let $(M,g,f)$ be complete with bounded Ricci curvature.  Assume
%$Rc_{f} \geq \lambda g$ with $\lambda > 0$.  Then there exists a
%compact smooth submanifold $K\subset\subset M$ with boundary and
%finite fundamental group such that each connected component of $M-K$
%is diffeomorphic to some $N\times\mathds{R}$, for $N$ some
%$(n-1)$-dimensional smooth compact manifold.
%\end{thm2}
%\begin{remark}
%In the above theorem our compact $(n-1)$-manifold $N$ may very well
%be different on different components of $M-K$.
%\end{remark}

Next we wish to use the above to understand the structure of
shrinking solitons under some simplified conditions.  We will prove
the following:

\begin{thm2}
Let $(M,g,f)$ be complete with $Rc+\nabla^{2}f=\lambda g$,
$\lambda>0$.  Assume $Rc\geq 0$ and that $f$ is either convex or
concave. Then $(M,g)$ is isometric to a finite quotient of
$E\times\mathds{R}^{k}$ where $E$ is a compact simply connected
Einstein manifold.
\end{thm2}
\begin{remark}
The point of the above is that the soliton structure on such an $M$
must be trivial. The noncompactness of $M$ must result purely from
an isometric $\mathds{R}^{n}$ factor, and on the compact component
$f$ behaves trivially.  This is not true in the case $\lambda=0$,
and there are nontrivial soliton structures on such manifolds (for
instance the cigar and Bryant solitons).
\end{remark}

To prove the above we introduce the notion of the $f$-Laplacian of a
function $u$.  The motivation comes directly from the standard
Laplace-Beltrami operator, which is defined as $\triangle =
\nabla^{*}\nabla$ with $\nabla^{*}$ the adjoint of the covariant
derivative with respect to the Riemannian volume form. Similarly we
define:

\begin{definition}
The $f$-Laplacian of a function $u$ is defined by
$\triangle_{f}u\equiv \nabla^{*f}\nabla u$, where the adjoint is
taken with respect to the $f$-measure $e^{-f}dv_{g}$.
\end{definition}

Under a positivity assumption on $Rc_{f}$ we have the following
estimate and Liouville type theorem:

\begin{thm3}
Let $(M,g,f)$ be complete with bounded Ricci curvature. Assume
$Rc_{f} \geq \lambda g$ with $\lambda > 0$.  Let
$u:M\rightarrow\mathds{R}$.  Then

1) If $\triangle_{f}u=0$ $\exists$ $\alpha>0$ such that if $|u|\leq
e^{\alpha d(x,p)^{2}}$ for  some $p\in M$, then $u=constant$.

2) If $\triangle_{f}u\geq 0$ and $u\leq C'$ for some
$C'\in\mathds{R}$, then $u=constant$.
\end{thm3}

The above situation is not typical if we weaken the geometric
constraint on $Rc_{f}$ a little to just $Rc_{f}\geq 0$.  Just for
instructional sake we show that not only is the above not true, but
that comparison estimates in general, in particular any Harnack type
estimate, must depend on $f$ itself and thus $Rc_{f}\geq 0$ is not a
sufficient condition to control many \textit{a priori} estimates of
$\triangle_{f}$:

\begin{thm4}
There exists $(M,g,f^{k})$ complete with
$u^{k}:M\rightarrow\mathds{R}$, $k\in\mathds{N}$ such that

1) $Rc+\nabla^{2}f^{k}\geq 0$ for any $k$.

2) $\triangle_{f^{k}}u^{k}=0$ on $M$ with $u^{k}>0$.

3) There exists $x,y\in M$ fixed such that $lim_{k}
\frac{u^{k}(x)}{u^{k}(y)} = \infty $.

\end{thm4}
\begin{remark}
The above example even has $Rm\geq 0$.  So even under a combination
of nonnegativity assumptions on $Rm$ and $Rc_{f}$ we do not have
complete \textit{a priori} control over solutions of the
$f$-Laplacian.
\end{remark}

\section{Proof of Theorems 1}

The key to understanding the metric measure spaces $(M,g,f)$ under
geometric assumptions on the $Rc_{f}$ is to control and understand
the behavior of $f$. On that note we begin by proving the following
estimate for $f$.

\begin{lemma1}
Let $(M,g,f)$ be complete with $|Rc|\leq C$ and $Rc+\nabla^{2}f\geq
\lambda g$ for $\lambda\in\mathds{R}$.  Let $\gamma:[0,L]\rightarrow
M$ be a weakly minimizing unit speed geodesic, $\gamma (0)=p\in M$.
Then $\nabla_{\dot{\gamma}} f(L)\geq\lambda L+a$ and $f(L)\geq
\frac{\lambda}{2}L^{2}+aL+b$, where $a=a(\lambda,n,C,f|_{B(p,2)})$
and $b=b(f|_{B(p,2)})$ and $B(p,2)$ is the geodesic ball of radius
$2$ centered at $p$.
\end{lemma1}

\begin{remark}
Note in the above that $a$ and $b$ do not depend on $L$, and hence
the estimate on $f$ actually holds for all $t\in [0,L]$
\end{remark}

\begin{proof}
First assume $L\geq 2$.  Let $E^{i}(p)$ be an orthonormal basis at
$p$ with $E^{n} = \dot{\gamma}$.  Define $E^{i}(t)$ as the parallel
transport of $E^{i}$ over $\gamma(t)$.  Let $h:[0,L]\rightarrow
\mathds{R}$ be Lipschitz with $h(0)=h(L)=0$ and
$Y^{i}(t)=h(t)E^{i}(t)$. Now for some $i$ let $\gamma_{s}(t):[0,L]
\rightarrow M$ be a 1-parameter family of curves with
$\gamma_{0}=\gamma$ and $\frac{d}{ds}\gamma_{s}=Y^{i}$ its
variation.

If $l(\gamma_{s})$ is defined as the length of $\gamma_{s}$ we have,
because $\gamma$ is a weakly minimizing geodesic, by the usual
second variation formula that

\[
0\leq \frac{d^{2}}{ds^{2}}l(\gamma_{s}) =
\int_{0}^{L}|\nabla_{\dot{\gamma}}Y^{i}|^{2} -
R(Y^{i},\dot{\gamma},Y^{i}, \dot{\gamma})dt
\]
\[
=\int_{0}^{L}(h')^{2}|E^{i}|^{2}+h^{2}
|\nabla_{\dot{\gamma}}E^{i}|^{2}
-h^{2}R(E^{i},\dot{\gamma},E^{i},\dot{\gamma})dt
\]
\[
=\int_{0}^{L}(h')^{2}-h^{2}R(E^{i},\dot{\gamma},E^{i},\dot{\gamma})dt
.
\]
Because this holds for each $i$ we can sum and use our assumption
that $Rc\geq -\nabla^{2}f+\lambda g$ to get

\[
\lambda\int_{0}^{L}h^{2}dt \leq (n-1)\int_{0}^{L}(h')^{2}dt
+\int_{0}^{L}h^{2}\nabla_{\dot{\gamma}\dot{\gamma}}^{2}f dt.
\]
Now we define $h$  by the formula
$$h(t) = \left\{\begin{array}{lr}
t &   0\leq t  \leq 1\\
1 &   1\leq t  \leq L-1\\
L-t & L-1\leq t\leq L
\end{array}\right.
 $$

Inserting yields
\[
\lambda(L-\frac{4}{3}) \leq (n-1)(2) + \int_{0}^{L}
\nabla_{\dot{\gamma}\dot{\gamma}}^{2}f dt
-\int_{0}^{L}(1-h^{2})\nabla_{\dot{\gamma}\dot{\gamma}}^{2}f dt.
\]
Since $\gamma$ is a unit speed geodesic and
$\nabla_{\dot{\gamma}\dot{\gamma}}^{2}f \geq \lambda-C$ we get

\[
\lambda L-2(\frac{2}{3}\lambda+(n-1)+\frac{2}{3}(C-\lambda))
+\nabla_{\dot{\gamma}} f(0) \leq \nabla_{\dot{\gamma}} f(L)
\]
or $\lambda L+a\leq \nabla_{\dot{\gamma}} f(L) \Rightarrow
\frac{\lambda}{2}L^{2}+a L +b \leq f(L)$.  Now this is for $L\geq
2$.  If we replace $a$ and $b$ by
$a'=a-2\lambda-sup_{B_{(p,2)}}|\nabla f|$ and
$b'=b-2\lambda-2a'-sup_{B_{(p,2)}}|f|$ then our inequality holds for
all $L$.
\end{proof}

The main application of this estimate is the following corollary,
which is a direct consequence of the remark following the lemma.
 It gives us a global lower bound on $f$ as a quadratic of the distance
function.

\begin{cor1}
Let $(M,g,f)$ be a complete manifold with $|Rc|\leq C$ and
$Rc+\nabla^{2}f\geq \lambda g$ for $\lambda\in\mathds{R}$, and let
$p\in M$. Then $\forall x\in M$ we have $f(x)\geq
\frac{\lambda}{2}d(x,p)^{2}+a d(x,p)+b$, where $d(x,p)$ is the
distance function to $p$ and $a,b$ depend only on the constants
$\lambda$,$n$,$C$ and $f|_{B(p,2)}$
\end{cor1}

The most important special case of the above is when $\lambda>0$.
The quadratic growth estimate on $f$ in this case immediately gives
us the following useful facts:

\begin{cor2}
Let $(M,g,f)$ be complete with bounded Ricci curvature and satisfy
$Rc+\nabla^{2}f\geq \lambda g$ with $\lambda>0$.  Then $f$ is
bounded below and proper.
\end{cor2}

We may now prove Theorem 1:

\begin{proof}[Proof of Theorem 1]
Using exponential coordinates at $p$ we have, since $Rc\geq -C$, by
the standard comparison that $\sqrt{detg}\leq
(sinh(\sqrt{C}r))^{n-1}\lesssim e^{(n-1)\sqrt{C}r}$.  Here by
definition it is understood that we say $s\lesssim t$ if $s\leq A t$
and $A$ is a constant depending only on the dimension and other
fixed variables, in this case just the dimension.  Now integrating
in the tangent space, where it is understood that $\sqrt{detg}(x)=0$
if $x$ is outside the segment domain of $p$, we have

\[
Vol_{f}(M)=\int_{S^{n-1}}\int_{0}^{\infty}e^{-f}\sqrt{detg} dr
ds_{n-1} \lesssim \int_{0}^{\infty}
e^{-\frac{\lambda}{2}r^{2}+(a+(n-1)\sqrt{C})r+b}dr<\infty.
\]

To see that the fundamental group is finite we lift to the universal
cover $\tilde{M}$.  Apply the above to see $\tilde{M}$ must also
have finite $f$-volume, because it too satisfies the geometric
conditions of the theorem. But this is impossible unless the order
of $\pi_{1}(M)$ is finite.
\end{proof}

%The second theorem is also an immediate consqeuence of Lemma 1:

%\begin{proof}[Proof of Theorem 2]
%First we will construct the compact subset $K$ and then we will
%prove that $K$ has finite fundamental group.  Fix $p\in M$ and
%recall the inequalities

%\begin{equation}\label{Ineq1}
%f(x)\geq \frac{\lambda}{2}d(x,p)^2+a d(x,p)+b
%\end{equation}

%\begin{equation}\label{Ineq2}
%\nabla_{\dot{\gamma}} f(x)\geq\lambda d(x,p)+a
%\end{equation}
%for some real numbers $a$ and $b$.  Corollary 2 tells us that $f$ is
%proper and bounded uniformly from below.  The second inequality
%above ensures that we may pick $r$ large enough so that on
%$M-B(p,r)$, $f$ has no critical points.

%Now since $f$ tends to infinity by the first inequality we can pick
%$c>0$ large enough so that $B(p,r)\subset f^{-1}(-\infty,c]$. We
%define $K\equiv f^{-1}(-\infty,c]$.  By the lower boundedness and
%properness of $f$ we see that $K$ is compact, and in particular so
%is $f^{-1}(c)$. For $c'\geq c$ we have, since $f$ has no critical
%points on $M-B(p,r)$, that the level sets $f^{-1}(c')$ are smooth
%compact $(n-1)$ dimensional manifolds, though possibly with several
%components. The gradient flow of $f$ then gives us a natural
%diffeomorphism between $f^{-1}(c)\times\mathds{R}$ and
%$f^{-1}(c,\infty)=K^{c}$, the complement of $K$.  To see that $K$
%has finite fundament group notice that $f^{-1}(c,\infty)$ retracts
%onto $f^{-1}(c)$, which is the boundary of $K$.  Thus
%$\pi_{1}(M)=\pi_{1}(K\cup f^{-1}(c,\infty)) = \pi_{1}(K\cup
%f^{-1}(c)) = \pi_{1}(K)$.

%\end{proof}

\section{Proofs of Theorems 2 and 3}
To prove Theorem 2 we begin by proving some results involving the
$f$-Laplace operator.

\begin{proof}[Proof of Theorem 3 (1)]
Let $x\in M$ be arbitrary.  Note by multiplying by $e^{-f}$ we get
\begin{equation}
\nabla^{i}(e^{-f}\nabla_{i}u)=0.
\end{equation}
Let $\phi: M\rightarrow\mathds{R}$ be a cutoff function with
$$\phi = \left\{\begin{array}{lr}
1 \textrm{ on }  B(x,1)\\
0\leq\phi\leq 1 \textrm{ on } B(x,1+r)-B(x,1)\\
0 \textrm{ on } M-B(x,1+r)
\end{array}\right.
 $$
where $r>0$ and $|\nabla\phi|\leq\frac{C}{r}$ for some $C$.
Multiplying the above by $\phi^{2}u$ and integrating we get
\[
-\int(2\phi u\nabla^{i}\phi\nabla_{i}u+\phi^{2}|\nabla u|^{2})
e^{-f}dv_{g} = 0
\]
\[
\int\phi^{2}|\nabla u|^{2} e^{-f}dv_{g} = -2\int(\phi
u\nabla^{i}\phi\nabla_{i}u) e^{-f}dv_{g}
\]
\[
\leq \int (\frac{1}{2}\phi^{2}|\nabla u|^{2}
+2u^{2}|\nabla\phi|^{2}) e^{-f}dv_{g}
\]
so that
\[
\int_{B(x,1)}|\nabla u|^{2} e^{-f}dv_{g}\leq 4\int_{M}u^{2}|\nabla
\phi|^{2} e^{-f}dv_{g}
\]
\[
\leq \frac{4C^{2}}{r^{2}}\int_{B_{1+r}-B_{1}}u^{2}e^{-f}dv_{g} \leq
\frac{4C^{2}}{r^{2}}\int_{M}u^{2}e^{-f}dv_{g}.
\]
But let $\alpha<\frac{\lambda}{4}$, and thus $u^{2}(x)\lesssim
e^{2\alpha d(x,p)^{2}}$. So in exponential coordinates we compute
\[
\int_{M}u^{2}e^{-f}dv_{g}\leq\int_{S^{n-1}}\int_{0}^{\infty}
e^{-(\frac{\lambda}{2}-2\alpha)r^{2}+a r+b}dr ds_{n-1} < \infty
\]
for some constants $a$ and $b$.  Thus we can tend
$r\rightarrow\infty$ to get
\[
\int_{B(x,1)}|\nabla u|^{2} e^{-f}dv_{g} = 0
\]
Since $x$ was arbitrary, $|\nabla u|=0$ and thus $u$=constant.
\end{proof}
\begin{proof}[Proof of Theorem 3 (2)]
This is much the same.  Since $u$ is bounded above we can assume, by
adding a constant, that $sup$ $u = 1$. Let $u^{+}(x)=max(u(x),0)$.
Let $x\in M$ such that $u(x)>0$ and $\phi$ as in the last part with
center $x$. Then our equation $\nabla^{i}(e^{-f}\nabla_{i}u)\geq 0$
gives

\[
-\int(2\phi u^{+}\nabla^{i}\phi\nabla_{i}u+\phi^{2}\nabla^{i}
u^{+}\nabla_{i} u) e^{-f}dv_{g} \geq 0
\]
so that
\[
\int_{M}\phi^{2}|\nabla u^{+}|^{2} e^{-f}dv_{g} \leq
\frac{4C^{2}}{r^{2}}\int_{M}(u^{+})^{2}e^{-f}dv_{g}
\]
But $u^{+}$ is bounded and $\int_{M}e^{-f}dv_{g}$ is finite.  So we
may limit out, using monotone convergence, to get $\int_{M}|\nabla
u^{+}|^{2} e^{-f}dv_{g}=0$.  So $u^{+}$ is constant.  Since
$u(x)>0$, $u$ is constant.
\end{proof}

\begin{proof}[Proof of Theorem 4]
Let $M=\mathds{R}\times S^{n}$ with the standard product metric. Let
$f^{k}(t,s)=k t$ for $k$ a constant, $t\in\mathds{R}$ and $s\in
S^{n}$.  Clearly $Rc+\nabla^{2}f\geq 0$.  Looking for a solution of
$\triangle_{f^{k}}u^{k}=0$ which is a function of only $t$ as well
we find $u^{k}(t)$ must satisfy the ode $u_{tt}-k u_{t}=0$.  So
$u^{k}(t,s)=e^{kt}$ is a solution of our equation.  Notice
$u^{k}>0$.  Let $x$ be on the $t=1$ slice and $y$ on the $t=0$
slice, then we see that $\frac{u^{k}(x)}{u^{k}(y)}=e^{k}$.  Let $k$
tend to $\infty$.
\end{proof}

\begin{remark}
The key point in the above is that $f$ is only well defined up to a
linear function on $M$, thus if $f$ itself does not growth faster
than linearly then we can add large linear terms to $f$ which will
have a significant impact on the solutions of the $f$-Laplacian.
\end{remark}

Now we apply the above to prove Theorem 3:

\begin{proof}[Proof of Theorem 3]
First we assume $f$ is convex.  Then $Rc\leq Rc+\nabla^{2}f =
\lambda g$.

The following computation is useful:
\[
\nabla^{i}R_{ij}+\nabla^{i}\nabla_{i}\nabla_{j}f=0
\]
\[
\frac{1}{2}\nabla_{j}R + \nabla_{j}(-R-n\lambda) +
Rc_{jk}\nabla^{k}f = 0
\]
\begin{equation}
\nabla_{i}R = 2Rc_{ij}\nabla^{j} f
\end{equation}
Now if we take the divergence of this we get
\begin{equation}
\triangle_{f}R = 2(\lambda R-|Rc|^{2})
\label{R_eq}
\end{equation}
A similar computation gives us
\begin{equation}
\triangle_{f}|Rc|^{2} = 2|\nabla Rc|^{2}+4(\lambda
|Rc|^{2}-R_{ijkl}R^{ik}R^{jl}) \label{Rc_eq}
\end{equation}

Now if $\partial_{i}$ is an eigenbasis for $Rc$ we write the rhs of
(\ref{R_eq}) as $(\lambda R-|Rc|^{2}) = \Sigma
R_{ii}(\lambda-R_{ii})\geq 0$ under our assumptions.  In particular
the scalar curvature is a bounded subsolution to $\triangle_{f}$,
and thus must be constant.  Plugging this in we see that $\Sigma
R_{ii}(\lambda-R_{ii})= 0$, which under our assumptions implies that
each term is zero and thus every eigenvalue of $Rc$ is either $0$ or
$\lambda$. By continuity the number of $0$ eigenvalues must be
constant, and thus $|Rc|=const$.  So writing $(\lambda
|Rc|^{2}-R_{ijkl}R^{ik}R^{jl})=\Sigma_{p,q}sec(\partial_{p},\partial_{q})
R_{pp}(\lambda-Rc_{qq})=0$ we then see that we must have $|\nabla
Rc|=0$ from (\ref{Rc_eq}).

By using deRham's Theorem we then have an isometric splitting of the
universal cover (which is a finite cover by the previous theorem)
into $E\times N$, where $E$ has Einstein constant $\lambda$ and $N$
is simply connected and Ricci flat.  By restricting $f$ to $N$ we
see $N$ has a soliton structure. We finally show that $N$ is
$\mathds{R}^{k}$.

We know by Ricci flatness that on $N$, $\nabla^{2}f=\lambda g$.  Now
by Lemma 1 we know $f$ always has a global minimum point, say $p\in
N$.  If $x\in N$ and $\gamma$ a geodesic connecting $x$ to $p$ we
see by integration that $\nabla_{\dot{\gamma}}f(x)=\lambda d(x,p)$
and $f(x)= \frac{\lambda}{2}d(x,p)^{2}+f(p)$ . Hence $f$ has a
unique nondegenerate minimum point at $p$.  Now we compute
\begin{equation}
R_{ijkq}\nabla^{q}f = (\nabla_{i}\nabla_{j}\nabla_{k}f -
\nabla_{j}\nabla_{i}\nabla_{k}f) = \nabla_{j}R_{ik}-\nabla_{i}R_{jk}
= 0. \label{Rm}
\end{equation}

In particular, because p is a nondegenerate critical point, for any
unit vector $X\in T_{p}M$ we can find $x\rightarrow p$ such that
$\frac{\nabla f(x)}{|\nabla f|}\rightarrow X$ (just use Taylor's
theorem in normal coordinates to see this). Dividing both sides of
(\ref{Rm}) by $|\nabla f|$,taking $X=\partial_{l}$ and limiting out
we get that $Rm(p)=0$. A final computation now gives us that
\[
\nabla_{\nabla f}|Rm|^{2} = \nabla^{p}f\nabla_{p}|Rm|^{2} =
2\nabla^{p}f R^{ijkl}\nabla_{p}R_{ijkl}
\]
\[
= -2\nabla^{p}f R^{ijkl} (\nabla_{i}R_{jpkl} + \nabla_{j}R_{pikl})
\]
\[
= -2R^{ijkl}(\nabla_{i}( \nabla^{p}f R_{jpkl})+\nabla_{j}(
\nabla^{p}f
R_{pikl})-R_{jpkl}\nabla_{i}\nabla^{p}f-R_{pikl}\nabla_{j}\nabla^{p}f)
\]
\[
= 4\lambda R^{ijkl}(R_{jikl}) = -4\lambda |Rm|^{2} \leq 0
\]

Our explicit formula for $f$ tells us that the negative gradient
flow from any $x\in N$ converges to $p$, and hence from the above
$|Rm|$ takes a maximum at $p$.  But we showed $Rm(p)=0$.  Hence
$Rm=0$.  Since $N$ is simply connected, $N$ is isometric to
$\mathds{R}^{k}$.

If we now instead assume that $f$ is concave the situation is more
simple.  We get $Rc\geq Rc+\nabla^{2}f = \lambda g$.  In particular
$M$ is compact and so $f$ must be a constant and so we are done.
\end{proof}
\section{Acknowledgements}
I would like to thank Gang Tian for guidance and tutelage during
this project.
%---------------------------------------------------------------------------------


\begin{thebibliography}{1}

  \bibitem{BE} D. Bakry and M. \'{E}mery, Diffusions
  Hypercontractives, in: {\em S\'{e}minaire de probabilit\'{e}s} XIX,
  1983/84, 177-206, Lecture Notes in Math. 1123, Springer, Berlin,
  1985.

  \bibitem{L}  J. Lott, Some geometric properties of the
  Bakry-\'{E}mery-Ricci tensor, {\em Comment. Math. Helv.} 78(2003)
  865-883.

  \bibitem{P} G. Perelman, Ricci flow with surgery on
  three-manifolds, {\em arXiv:math/0303109}, 03/2003.

  \bibitem{MT} J. Morgan and G. Tian, Ricci Flow and the
  Poincar\'{e} Conjecture, {\em arXiv:math/0607607}, 07/2006.

  \end{thebibliography}
\end{document}